\documentclass{amsart}
\usepackage{graphicx} 
\usepackage{amsmath,amssymb, amsthm,mathrsfs}
\usepackage{amsmath}
\usepackage[all]{xy}
\usepackage{hyperref}

\usepackage[a4paper,top=2cm,bottom=2cm,left=3cm,right=3cm,marginparwidth=1.75cm]{geometry}

\newtheorem{theorem}{Theorem}[section]

\theoremstyle{definition}
\newtheorem{definition}[theorem]{Definition}

\newtheorem{proposition-definition}[theorem]{Proposition-Definition}

\begin{document}

\title{A short proof for the acyclicity of oriented exchange graphs of cluster algebras}

\author{Shuhao Deng and Changjian Fu}

\address{Department of Mathematics\\SiChuan University\\610064 Chengdu\\P.R.China}
\email{878915687@qq.com(Deng)\\ changjianfu@scu.edu.cn(Fu)}

\begin{abstract}
The statement in the title was proved in \cite{Cao23} by introducing dominant sets of seeds, which are analogs of torsion classes in representation theory. In this note, we observe a short proof by the existence of consistent cluster scattering diagrams.
\end{abstract}
\maketitle
\section{Recollection}
\subsection{Oriented exchange graph of a cluster algebra}
We refer to \cite{FZ02, FZ07} for the unexplained terminology in cluster algebras. Let $r$ be a positive integer and $\mathbf{\Sigma}=\{\Sigma_t=(\mathbf{x}_t,B_t)\}_{t\in \mathbb{T}_r}$ a cluster pattern of rank $r$ with trivial coefficients, where $\mathbb{T}_r$ is the $r$-regular tree. We fix a vertex $t_0\in \mathbb{T}_r$ as the root vertex and $\Sigma_{t_0}$ as the initial seed. Denote by $B:=B_{t_0}$.

For each vertex $t\in \mathbb{T}_r$, denote by $C_t^{t_0}=(\mathbf{c}_{1;t}^{t_0},\dots, \mathbf{c}_{r;t}^{t_0})\in M_r(\mathbb{Z})$ the {\it $C$-matrix} at $t$ with respect to the root vertex $t_0$. The column vectors $\mathbf{c}_{1;t}^{t_0},\dots, \mathbf{c}_{n;t}^{t_0}$ are {\it $c$-vectors} at $t$ with respect to $t_0$. 
Denote by $G_t^{t_0}=(\mathbf{g}_{1;t}^{t_0},\dots, \mathbf{g}_{r;t}^{t_0})\in M_r(\mathbb{Z})$ the {\it $G$-matrix} at $t$ with respect to $t_0$. Its column vectors $\mathbf{g}_{1;t}^{t_0},\dots, \mathbf{g}_{r;t}^{t_0}$ are {\it $g$-vectors} at $t$ with respect to $t_0$. Both $C_t$ and $G_t$ are invertible over $\mathbb{Z}$ and we have the so-called {\it tropical duality} \cite{NZ12}
\[
(G_t^{t_0})^{\operatorname{tr}}DC_t^{t_0}=D,
\]
where $D$ is a skew-symmetrizer of the exchange matrix $B_{t_0}$.
Recall that a non zero integer vector $\alpha=(a_1,\dots, a_r)'\in\mathbb{Z}^r$ is {\it positive} (resp. {\it negative}) if $a_i\geq 0$ (resp. $a_i\leq 0$) for all $1\leq i\leq r$. In either cases, we say that $\alpha$ is {\it sign-coherent}. It has been proved by \cite{GHKK18} that $\mathbf{c}_{i;t}^{t_0}$ is sign-coherent for each $t\in \mathbb{T}_r$ and $1\leq i\leq r$. As a consequence, the following is well defined.

\begin{definition}
    The mutation $\mu_k(\Sigma_t)$ of $\Sigma_t$ in direction $k$ is a {\it green mutation} (resp. {\it red}) if $\mathbf{c}_{k;t}^{t_0}$ is positive (resp. negative).
\end{definition}
Denote by $[\Sigma_t]$ the equivalence class of $\Sigma_t$ up to simultaneously relabeling the $r$-tuple $\mathbf{x}_t$ and the rows and columns of $B_t$. We refer to $[\Sigma_t]$ as an {\it unlabeled seed}. It is known that the mutation $\mu_k(\Sigma_t)$ induces a mutation $\mu_{x_{k;t}}([\Sigma_t])=[\mu_k(\Sigma_t)]$ of unlabeled seeds. Furthermore, each mutation is either green or red.
\begin{definition}
    The {\it oriented exchange graph} or {\it exchange quiver} $\vec{\mathcal{H}}_{t_0}(\mathbf{\Sigma})$ of the cluster pattern $\mathbf{\Sigma}$ has vertices consisting of unlabeled seeds of $\mathbf{\Sigma}$ and there is an arrow from $[\Sigma_t]$ to $[\Sigma_{t'}]$ if and only if $[\Sigma_{t'}]$ is a green mutation of $[\Sigma_t]$.
\end{definition}
    By definition, $[\Sigma_{t_0}]$ is a source vertex in $\vec{\mathcal{H}_{t_0}}(\mathbf{\Sigma})$.
    The following has been proved by \cite[Theorem 1.6(i)]{Cao23} by introducing dominant sets of seeds, which are analogs of torsion classes in representation theory.
    \begin{theorem}\label{thm:acyclic}
        The oriented exchange graph $\vec{\mathcal{H}}_{t_0}(\mathbf{\Sigma})$ is acyclic.
    \end{theorem}
The aim of this note is to show that Theorem \ref{thm:acyclic} follows from the existence of cluster scattering diagrams \cite{GHKK18}.

\subsection{Cluster scattering diagram}
Cluster scattering diagram was introduced in \cite{GHKK18}, which is a powerful tool to study cluster algebras. Here we mainly follow the notations in \cite{Na23}.
\subsubsection{The structure group}
Fix a decomposition $B=\Delta\Omega$, where $\Delta=\operatorname{diag}\{\delta_1,\dots, \delta_r\}$ with positive integers $\delta_1,\dots, \delta_r$ and $\Omega=(w_{ij})$ is an $r\times r$ skew-symmetric rational matrix. Let $N$ be a rank $r$ lattice with basis $e_1,\dots,e_r$ and $\{-,-\}:N\times N\to \mathbb{Q}$ a skew-symmetric $\mathbb{Z}$-bilinear form given by $\{e_i,e_j\}=w_{ij}$. Let $N^o$ be the sublattice of $N$ generated by $\delta_1e_1,\dots, \delta_re_r$. Denote by $M=\operatorname{Hom}(N,\mathbb{Z})$ and $M^\circ=\operatorname{Hom}(N^\circ,\mathbb{Z})$ the duals of $N$ and $N^\circ$ respectively. Then $\Gamma=(N,\{-,-\}, N^\circ,\Delta, M,M^\circ)$ is called a {\it data} of $B$ and $\mathfrak{s}=(e_1,\dots, e_r)$ a {\it seed} of $\Gamma$.

Let $N^+=\{\sum_{i=1}^ra_ie_i~|~a_i\in \mathbb{Z}_{\geq 0}, \sum_{i=1}^ra_i>0\}$ be the set of positive elements in $N$ and $N_{\text{pr}}^+$ the set of primitive elements of $N^+$. The {\it degree function} $\deg: N^+\to \mathbb{Z}_{>0}$ is defined by $\deg(\sum_{i=1}^ra_ie_i)=\sum_{i=1}^ra_i$.

Let $\mathfrak{g}=\bigoplus_{n\in N^+}\mathbb{C}X_n$ be the $N^+$-graded $\mathbb{C}$-vector space endowed with a basis $\{X_n\mid n\in N^+\}$. It is an $N^+$-graded Lie algebra if we define the Lie bracket by $[X_n,X_{n'}]=\{n,n'\}X_{n+n'}$. For any positive integer $l$, \[\mathfrak{g}^{>l}:=\bigoplus_{n\in N^+, \deg{n}>l}\mathbb{C}X_n\] is an ideal of $\mathfrak{g}$ and the quotient Lie algebra $\mathfrak{g}^{\leq l}:=\mathfrak{g}/\mathfrak{g}^{>l}$ is nilpotent. For any $l'\geq l$, we have a canonical projection $\pi_{l',l}:\mathfrak{g}^{\leq l'}\to \mathfrak{g}^{\leq l}$, and hence we can form the completion of $\mathfrak{g}$:
\[
\hat{\mathfrak{g}}=\lim_{\leftarrow}\mathfrak{g}^{\leq l}.
\]
Denote by $\pi_l:\hat{\mathfrak{g}}\to \mathfrak{g}^{\leq l}$ the canonical projection for any positive integer $l$. Let $G^{\leq l}=\operatorname{exp}(\mathfrak{g}^{\leq l})$ be the exponential group of $\mathfrak{g}^{\leq l}$. Similarly, we have $G:=\lim\limits_{\leftarrow}G^{\leq l}$ and the canonical projection $\pi_l:G\to G^{\leq l}$. The group $G$ is called the {\it structure group} of the data $\Gamma$ and $\mathfrak{s}$.


For $n\in N_{\text{pr}}^+$, let $\mathfrak{g}_n^{\parallel}:=\bigoplus_{j>0}\mathbb{C}X_{jn}$, which is an abelian Lie subalgebra of $\mathfrak{g}$. Denote by $G_n^{\parallel}$ the abelian subgroup of $G$ corresponding to the completion of $\mathfrak{g}_n^{\parallel}$. The group $G_n^{\parallel}$ is called the {\it parallel subgroup of $n$}.

\subsubsection{Cluster scattering diagram}

Denote by $M_\mathbb{R}=M\otimes_\mathbb{Z} \mathbb{R}$ and $\langle-,-\rangle: N\times M_\mathbb{R}\to \mathbb{R}$ the canonical pairing.

A {\it wall} for the seed $\mathfrak{s}$ is a triplet $\mathbf{w}=(\mathfrak{d},g)_n$, where $n\in N^+_{\text{pr}}$, $\mathfrak{d}$ is a convex rational polyhedral cone spanning $n^\perp=\{z\in M_\mathbb{R}\mid \langle n,z\rangle=0\}$ and  $g\in G_n^{\parallel}$.  We call $n,\mathfrak{d},g$ the {\it normal vector}, the {\it support}, the {\it wall element} of $\mathbf{w}$.
The support $\mathfrak{d}$ is also referred to as a wall.
For a wall $\mathbf{w}=(\mathfrak{d},g)_n$, the open half space $\{m\in M_\mathbb{R}\mid \langle n,m\rangle>0\}$ is the {\it green side} of $\mathbf{w}$, and $\{m\in M_\mathbb{R}\mid \langle n,m\rangle<0\}$ is the {\it red side} of $\mathbf{w}$.
\begin{definition}
    A {\it scattering diagram $\mathfrak{D}$} for $\mathfrak{s}$ is a collection of walls $\{\mathbf{w}_\lambda=(\mathfrak{d}_\lambda,g_\lambda)_{n_\lambda}\}_{\lambda \in \Lambda}$, where $\Lambda$ is a finite or countably infinite index set, satisfying the finiteness condition: For each positive integer $l$, there are only finitely many walls such that $\pi_l(g_\lambda)\neq \operatorname{id}$.
\end{definition}
Each connected compartment of $M_{\mathbb{R}}\backslash\mathfrak{D}$ is called a {\it chamber} of $\mathfrak{D}$.

Let $\mathfrak{D}=\{\mathbf{w}_\lambda=(\mathfrak{d}_\lambda,g_\lambda)_{n_\lambda}\}$ be a scattering diagram of $\mathfrak{s}$. A smooth curve $\gamma:[0,1]\to M_\mathbb{R}$ is admissible if the following conditions are satisfied:
\begin{itemize}
    \item $\gamma(0)$ and $\gamma(1)$ are not in any wall $\mathfrak{d}_\lambda$ of $\mathfrak{D}$;
    \item $\gamma$ crosses walls  of $\mathfrak{D}$ transversely;
    \item  $\gamma$ does not cross the boundaries of walls or any intersection of walls spanning two
distinct hyperplanes.
\end{itemize}

Let $\gamma$ be an admissible curve of $\mathfrak{D}$. For simplicity, we assume that $\gamma$ crosses finitely many walls of $\mathfrak{D}$ in the following order:
\[
\mathbf{w}_1=(\mathfrak{d}_1,g_1)_{n_1},\dots,\mathbf{w}_s=(\mathfrak{d}_s,g_s)_{n_s}.
\]
The {\it path-ordered product} of $\gamma$ is defined as 
\[
\mathfrak{p}_{\gamma,\mathfrak{D}}:=g_s^{\epsilon_s}\cdots g_1^{\epsilon_1}\in G,
\]
where $\epsilon_i=1$ (resp. $-1$) if $\gamma$ crosses $\mathfrak{d}_i$ from its green side to red side (resp. from its red side to its green side). In general, one can define the path-ordered product by limits.
\begin{definition}
    A scattering diagram $\mathfrak{D}$ is {\it consistent} if the path-ordered product $\mathfrak{p}_{\gamma,\mathfrak{D}}=\operatorname{id}$ for any admissible loop $\gamma$ of $\mathfrak{D}$.
\end{definition}

We say that $z\in M_{\mathbb{R}}$ is {\it general} if there is at most one rational hyperplane $n^\perp(n\in N,n\neq 0)$ such that $z\in n^\perp$.

For any $n\in N^+$, define
\[
\Psi[n]:=\operatorname{exp}(\sum_{j=1}^\infty\frac{(-1)^{j+1}}{j^2}X_{jn})\in G,
\]
which is called the {\it dilogarithm element} of $n$. Note that if $n\in N_{\text{pr}}^+$, then $\Psi[n]\in G_n^{\parallel}$. The scattering diagram in the following theorem is the {\it cluster scattering diagram} associated with the seed $\mathfrak{s}$.
\begin{theorem}[{\cite[Theorem 1.12, 1.13]{GHKK18}}]\label{thm:consistentdiagram}
There exists a unique consistent scattering diagram $\mathfrak{D}_0(B)$ (up to equivalence) for the seed $\mathfrak{s}$, satisfying the following conditions:
\begin{itemize}
\item For any $i\in\{1,\dots,r\}$, $(e_i^{\perp},\Psi[e_i]^{\delta_i})_{e_i}$ is a wall of $\mathfrak{D}_0(B)$.
\item For any other walls $(\mathfrak{d},g)_n$ of $\mathfrak{D}_0(B)$, $\{-,n\}\notin \mathfrak{d}$ holds. (if we take $(f_1,...,f_r)$ in the following section as the basis of $M_\mathbb{R}$, then $\{-,n\}$ is $Bn$.)
\end{itemize}
Moreover, we can construct $\mathfrak{D}_0(B)$ such that the wall element of any wall $\mathbf{w}=(\mathfrak{d},g)_n$ has the following form $g=\Psi[tn]^{s\delta(tn)}$, where $s,t\in \mathbb{Z}_{>0}$ and $\delta(tn)$ is the smallest positive rational number such that $\delta(tn)tn\in N^\circ$. 
\end{theorem}
In the following, for a convex cone $W$, we denote by $\operatorname{Int}(W)$ its interior and by $\partial W$ its boundary. 

\section{Proof of the statement}

Assume that $\vec{\mathcal{H}}_{t_0}(\mathbf{\Sigma})$ is not acyclic. In particular, there is an oriented cycle in $\vec{\mathcal{H}}_{t_0}(\mathbf{\Sigma})$, which yields a path in $\mathbb{T}_r$:
\begin{equation}\label{p:path}
  \xymatrix{t_1\ar@{-}[r]^{k_1}&t_2\ar@{-}[r]^{k_2}&\cdots\ar@{-}[r]^{k_{s-1}}&t_s\ar@{-}[r]^{k_s}&t_{s+1}}  
\end{equation}
such that $[\mathbf{x}_{t_1}]=[\mathbf{x}_{t_{s+1}}]$ and for $1\leq i\leq s$, $\mu_{k_i}(\Sigma_{t_i})$ is a green mutation of $\Sigma_{t_i}$. In particular, $\mathbf{c}_{k_i;t_i}^{t_0}$ is positive and $k_i\neq k_{i+1}$.

Let $e_1^\ast,\dots, e_r^\ast\in M$ be the dual basis of $e_1,\dots,e_r$. Denote by $f_i=\delta_i^{-1}e_i^\ast$ for $1\leq i\leq r$. It is clear that $f_1,\dots, f_r$ is an $\mathbb{R}$-basis of $M_\mathbb{R}$. We identify $M_\mathbb{R}$ with $\mathbb{R}^r$ by identifying $f_1,\dots, f_r$ with the standard basis of $\mathbb{R}^r$.

Recall that for each $t\in \mathbb{T}_r$, we have a $G$-matrix $G_t^{t_0}=(\mathbf{g}_{1;t}^{t_0},\dots, \mathbf{g}_{r;t}^{t_0})$. We denote by $\sigma(G_t^{t_0})$ the convex cone generated by the $g$-vectors at $t$, and by $\sigma_i(G_t^{t_0})$ the convex cone generated by the $g$-vectors at $t$ except its $i$th $g$-vector $\mathbf{g}_{i;t}^{t_0}$.

Let $\mathfrak{D}_0(B)$ be the consistent cluster scattering diagram associated with the seed $\mathfrak{s}$ in Theorem \ref{thm:consistentdiagram}. It has been proved by \cite[Lemma 2.10]{GHKK18} that $\operatorname{Int}(\sigma (G_t^{t_0}))$ is a chamber of $\mathfrak{D}_0(B)$ for any vertex $t\in \mathbb{T}_r$.

By the tropical duality, we know that for $1\leq i\leq s$,
\begin{itemize}
    \item $\sigma(G_{t_i}^{t_0})\cap \sigma(G_{t_{i+1}}^{t_0})=\sigma_{k_i}(G_{t_i}^{t_0})$,
    \item $n_i:=(e_1,\dots, e_r)\mathbf{c}_{k_i;t_i}^{t_0}$ is the normal vector of $\sigma_{k_i}(G_{t_i}^{t_0})$,
    \item $\operatorname{Int}(\sigma(G_{t_i}^{t_0}))$ (resp. $\operatorname{Int}(\sigma(G_{t_{i+1}}^{t_0}))$) is in the green (resp. red) side of $\sigma_{k_i}(G_{t_i}^{t_0})$.
\end{itemize}

According to the path \eqref{p:path}, we may construct an admissible loop $\gamma$ with an endpoint in $\operatorname{Int}(\sigma(G_{t_1}^{t_0}))$ and passing through $\operatorname{Int}(\sigma(G_{t_2}^{t_0})),\dots, \operatorname{Int}(\sigma(G_{t_s}^{t_0}))$ in sequence. The loop crosses general points $z_1\in \operatorname{Int}(\sigma_{k_1}(G_{t_1}^{t_0})),\dots, z_s\in \operatorname{Int}(\sigma_{k_s}(G_{t_s}^{t_0}))$ in order.

Note that a general point $z$ may lie in infinitely many walls with the same normal vector, but for any general point $z\in\sigma_{k_i}(G_{t_i}^{t_0})$, it lies on exactly one wall whose wall element is $\Psi[n_i]^{\delta(n_i)}$ \cite[Proposition 6.24]{Na23}. According to the definition of the path-ordered product of $\gamma$ and Theorem \ref{thm:consistentdiagram}, we have
\begin{equation*}
    \operatorname{id}=\mathfrak{p}_{\gamma,\mathfrak{D}_0(B)}=\Psi[n_s]^{\delta(n_s)}\Psi[n_{s-1}]^{\delta(n_{s-1})}...\Psi[n_1]^{\delta(n_1)}.
\end{equation*}
Let $l=\min\{\deg(n_i)\mid i=1,\dots, s\}$. By applying the canonical projection $\pi_l:G\to G^{\leq l}$, we obtain
\[
\operatorname{id}_{G^{\leq l}}=\pi_l(\operatorname{exp}(\sum_{i,\deg(n_i)=l}\delta(n_i)X_{n_i})).
\]
However, $\pi_l(\operatorname{exp}(\sum_{i,\deg(n_i)=l}\delta(n_i)X_{n_i}))\neq \operatorname{id}_{G^{\leq l}}$, a contradiction. This completes the proof.

\noindent{\bf Acknowledgment.} The authors are grateful to Prof. Peigen Cao for his interest and careful reading of this draft.

\end{document}